\documentclass[12pt, reqno]{amsart}
\usepackage{amsmath, amsthm, amscd, amsfonts, amssymb, graphicx, color}
\usepackage[bookmarksnumbered, colorlinks, plainpages]{hyperref}
\hypersetup{colorlinks=true,linkcolor=red, anchorcolor=green, citecolor=cyan, urlcolor=red, filecolor=magenta, pdftoolbar=true}
\newtheorem{theorem}{Theorem}[section]
\newtheorem{lemma}[theorem]{Lemma}
\newtheorem{proposition}[theorem]{Proposition}

\theoremstyle{definition}
\newtheorem{definition}[theorem]{Definition}
\newtheorem{example}[theorem]{Example}

\theoremstyle{remark}

\usepackage{mathtools}

\begin{document}

\setcounter{page}{1}

\title{  semi-order continuous operators on vector spaces}
\author[K. Haghnejad Azar]{Kazem Haghnejad Azar}
\author[M. Matin]{Mina Matin}
\author[R. Alavizadeh]{Razi Alavizadeh}
\address{Department  of  Mathematics  and  Applications, Faculty of  Sciences, University of Mohaghegh Ardabili, Ardabil, Iran.}
\email{\textcolor[rgb]{0.00,0.00,0.84}{haghnejad@uma.ac.ir}}
\email{\textcolor[rgb]{0.00,0.00,0.84}{minamatin@uma.ac.ir}}
\email{\textcolor[rgb]{0.00,0.00,0.84}{ralavizadeh@uma.ac.ir}}
\keywords{ ordered vector space, pre-Riesz space, $\tilde o$-convergence, $\tilde o$-continuous, $o_{W_1}$-convergence, $o_{W_1}$-continuous.
\newline \indent $^{*}$Corresponding author}
\begin{abstract}
In this manuscript, we will study both $\tilde{o}$-convergence in (partially) ordered vector spaces and  a kind of convergence in a vector space $V$.
A vector space $V$ is called semi-order vector space (in short semi-order space), if there exist an ordered vector space $W$ and an operator $T$ from $V$ into $W$. In this way, we say that $V$ is semi-order space with respect to $\{W, T\}$.
	 A net $\{x_\alpha\}\subseteq V$ is said to be ${\{W,T\}}$-order convergent  to a vector $x\in V$ (in short we write $x_\alpha\xrightarrow {\{W, T\}}x$), whenever there exists a net $\{y_\beta\}$ in $W$  satisfying $y_\beta \downarrow 0$ in $W$ and for each $\beta$, there exists $\alpha_0$ such that $\pm (Tx_\alpha -Tx) \leq y_\beta$ whenever $\alpha \geq \alpha_0$. In this manuscript, we study and investigate some properties of  $\{W,T\}$-convergent nets and its relationships with other order convergence in partially ordered vector spaces.
Assume that $V_1$ and $V_2$ are semi-order spaces with respect to 	 $\{{W_1}, T_1\}$ and $\{W_2, T_2\}$, respectively. An operator $S$ from $V_1$ into $V_2$ is called semi-order continuous, if  $x_\alpha\xrightarrow {\{{W_1}, T_1\}}x$ implies $Sx_\alpha\xrightarrow {\{W_2, T_2\}}Sx$ whenever $\{x_\alpha\}\subseteq V_1$. We study some properties of this new
 classification of  operators.

\end{abstract}
\maketitle

\section{Introduction}
Let $W$ be a real vector space and $K$ be a cone in $W$, that is, $K$ is a wedge ($x,y \in K$ and $\lambda,\mu \geq 0$ imply $\lambda x + \mu y \in K$) and $K \cap (-K) = \{0\}$. In $W$ a partial order is defined by $x\leq y $ whenever $y-x \in K$. The space $(W, K)$ (or, loosely $W$) is then called a (partially) ordered vector space. A subspace $M \subseteq W$ is majorizing in $W$ if for every $x\in W$ there exists some $m\in M$ with $x\leq m$ (or, equivalently, if for each $x\in W$ there exists some $m\in M$ with $m\leq x$). A subspace $M\subseteq W$ is called directed if for every $x,y \in M$ there is an element $z\in M$ such that $x\leq z$ and $y\leq z$. An ordered vector space $W$ is directed if and only if $W_+ $ is generating in $W$, that is, $W = W_+ - W_+$.
An ordered vector space $W$ is called Archimedean if for every $x,y \in W$ with $nx\leq y$ for every $n\in \mathbb{N}$ one has $x\leq 0$. 
The ordered vector space $W$ has the Riesz decomposition property (RDP) if for every $x_1 , x_2 , z \in K$ with $z\leq x_1 + x_2$ there exist $z_1 , z_2 \in K$ such that $ z = z_1 +z_2$ with $z_1 \leq x_1$ and $z_2 \leq x_2$. We call a linear subspace $M$ of an ordered vector space $W$ order dense in $W$ if for every $x\in W$ we have 
\begin{equation*}
x = \inf \{z\in M : x \leq z\},
\end{equation*}
that is, the greatest lower bound of set $\{z\in M : x \leq z\}$ exists in $W$ and equals $x$, see page $360$ of \cite{1a}. Clearly, if $M$ is order dense in $W$, then $M$ is majorizing in $W$. Denote for a subset $M$ of $W$, the set of all upper bounds by $M^u = \{x\in W : x\geq m \ for \ all\ m \in M\}$.
   A subset $M$ of an ordered vector space $W$ is called solid if for every $x\in W$ and $y\in M$ the relation $\{\pm y\}^u \subseteq \{\pm x\}^u$ implies $x\in M$. A solid subspace $M$ of $W$ is called an ideal.  The elements $x,y \in W$ are called disjoint, in symbols $x \perp y$, if $\{\pm (x+y)\}^u = \{\pm (x-y)\}^u$. The disjoint complement of a subset $M\subseteq W$ is $M^d = \{x\in W \mid \forall y\in M: x \perp y\}$. A sequence $\{x_n\}\subseteq W$ is said to be disjoint, if for every $n \neq m$, $x_n \perp x_m$.  A linear subspace $M$ of an ordered vector space $W$ is called a band in $W$ if $M = M^{dd} $.\\
 Recall that a linear map $i : W_1 \rightarrow {{W_2}}$ between two ordered vector spaces is said to be bipositive if for every $x\in W_1$ one has $i(x)\geq 0$ if and only if $x\geq 0$. A partially ordered vector space $(W_1,K)$ is called pre-Riesz space if for every $x,y,z \in W_1$ the inclusion $\{x+y,x+z\}^u \subseteq \{y,z\}^u$ implies $x\in K$. Clearly, each vector lattice is pre-Riesz space, since the inclusion in definition of pre-Riesz space reduces to inequality $(x+y)\vee (x+z)\geq y \vee z$, so $ x+(y \vee z) \geq y \vee z$, which implies $x\geq 0$. By Theorem 4.3 of \cite{3},  ordered vector space $W_1$ is a pre-Riesz space if and only if there exist a vector lattice ${{W_2}}$ and a bipositive linear map $i: W_1\rightarrow {{W_2}}$ such that $i(W_1)$ is order dense in ${{W_2}}$. The pair $({{W_2}},i)$ (or, loosely ${{W_2}}$) is then called a vector lattice cover of $W_1$. The theory of pre-Riesz spaces and their vector lattice covers is due to van Haandel, see \cite{5}.\\
A net $\{x_\alpha\} $ in ordered vector space $W$ is said to be decreasing (in symbols, $x_\alpha \downarrow$), whenever $\alpha \geq \beta $ implies $x_\alpha \leq x_\beta$. For $x\in W$ the notation $x_\alpha \downarrow x $ means that $x_\alpha \downarrow $ and $\inf_\alpha \{x_\alpha\} = x$ both hold. The meanings of $x_\alpha \uparrow x$ are analogus.  We say that a net $\{x_\alpha\}\subseteq W$, $(o)$-converges to $x\in W$ (in symbols, $x_\alpha \xrightarrow{o}x)$, if there is a net $\{y_\alpha\}\subseteq W$ such that $y_\alpha \downarrow 0$ and for all $\alpha$ one has $ \pm(x_\alpha - x) \leq y_\alpha$. Let $W_1$ and $W_2$ be two ordered vector spaces. An operator $T:W_1\rightarrow W_2$ is said to be order continuous, if every net $\{x_\alpha\}\subseteq W_1$ with $x_\alpha \xrightarrow{o}0 $ implies $Tx_\alpha \xrightarrow{o}0$ in $W_2$. The collection of all order continuous operators between two ordered vector spaces $W_1$ and $W_2$, showed by $L_{oc}(W_1,W_2)$. From \cite{2} recall that the directed part of $L_{oc}(W_1,W_2)$ has been showed by $L_{oc}^\diamond (W_1,W_2) = L_{oc}(W_1,W_2)_+ - L_{oc}(W_1,W_2)_+$ where $W_1$ and $W_2$ are ordered vector spaces and $W_1$ is directed.  \\
Let $(W,K)$ be an ordered vector soace. For two elements $y,z \in K$ with $ y \leq z$ denote the according order interval by $ [y,z] = \{x\in W: y\leq x \leq z\}$. A set $M\subset W$ is called order bounded if there are $y,z \in W$ such that $M\subseteq [y,z]$.\\

\section{ $\tilde o$-continuous operators}
	Let $W$ be an ordered vector space. A net $\{x_\alpha\}\subseteq W$ is said to be $\tilde{o}$-convergent to $x\in W$ (in symbol, $x_\alpha \xrightarrow{\tilde{o}} x $) if  there exists a net $\{y_\beta\}\subseteq W$, possibly over a different index set, such that $y_\beta \downarrow 0$ in $W$ and for every $\beta$ there exists $\alpha_0$ such that $\pm(x_\alpha - x)\leq y_\beta$, whenever $\alpha \geq \alpha_0$.\\	
It is clear that for each net $\{x_\alpha\}\subseteq W$, $x_\alpha \xrightarrow{o}x$ implies $x_\alpha \xrightarrow{\tilde{o}}x$, but by Example 1.4 of \cite{11}, the converse, in general not holds.

\begin{lemma}\label{eli}
	Let $W$ be an ordered vector space and $\{x_\alpha\}\subseteq W$. Then we have the following assertions.
\end{lemma}
\begin{enumerate}
	\item $x_\alpha \xrightarrow{\tilde o}x$ iff $(x_\alpha - x) \xrightarrow{\tilde o}0$.
	\item If $0\leq x_\alpha \xrightarrow{\tilde o}x$, then $0 \leq x$.
	\item If for each $\alpha$, $ x_\alpha \leq y$ and $x_\alpha \xrightarrow{\tilde o}x$, then $ x \leq y$.
	\item If $x_\alpha\xrightarrow{\tilde o}x$ and $x_\alpha \xrightarrow{\tilde o}y$, then $x = y$.
\item If $x_\alpha \xrightarrow{\tilde o}x$ and $y_\alpha \xrightarrow{\tilde o}y$, then $ \lambda x_\alpha + \mu y_\alpha \xrightarrow{\tilde o} \lambda x + \mu y$ for all $ \lambda , \mu  \in\mathbb{R}$.
\item If $x_\alpha \xrightarrow{\tilde o}x$, $z_\alpha \xrightarrow{\tilde o}z$ and $ x_\alpha \leq z_\alpha$ for all $\alpha$, then $ x \leq z$.
\end{enumerate}
\begin{proof}
	\begin{enumerate}
	\item By definition it is established.
	\item	Since $x_\alpha \xrightarrow{\tilde o}x$, therefore there exists net $\{y_\beta\}\subseteq W$ such that $y_\beta \downarrow 0$ and for every $\beta $ there exists $\alpha_0$ such that $ \pm (x_\alpha - x) \leq y_\beta$ whenever $ \alpha \geq \alpha_0$. Since  $x_\alpha - x \leq y_\beta$, then $0 \leq x_\alpha \leq x + y_\beta$. Therefore $0\leq y_\beta + x$ and $0 \leq \inf_\beta \{y_\beta\} + x = x$.
	\item By assumption there exists a net $\{y_\beta\}\subseteq W$ such that $y_\beta \downarrow 0$ and for every $\beta $ there exists $\alpha_0$ such that $ \pm (x_\alpha - x) \leq y_\beta$ whenever $ \alpha \geq \alpha_0$. We have $ x = x - x_\alpha + x_\alpha \leq y_\beta +y$ whenever $\alpha \geq \alpha_0$ and therefore $ x \leq \inf_\beta \{y_\beta\} + y$, and follows $ x \leq y$.
	\item We have $x - y = x - x_\alpha + x_\alpha - y \leq y_\beta + z_\gamma $ where $y_\beta \downarrow 0 , z_\gamma \downarrow 0$, therefore $x\leq y$. By the same way $ y \leq x$ and so $ x = y$.
	\item We know that if $y_\beta \downarrow 0$ and $x\leq y$, then for every $\lambda \geq 0$, $\lambda y_\beta \downarrow 0$ and $\lambda x \leq \lambda y$. Note that if $\lambda < 0$, then $\lambda (\pm x_\alpha) = -\lambda (\pm x_\alpha)$. Therefore the proof holds.
	\item According to the relationship $x - z = x - x_\alpha + x_\alpha - z \leq x - x_\alpha + z_\alpha - z$, the proof is obvious.\qedhere
	\end{enumerate}
\end{proof}
\begin{definition}
	Let $W_1$ and $W_2$ be two ordered vector spaces. An operator $T:W_1\rightarrow W_2$ is said to be
	\begin{enumerate}
		\item  $\tilde{o}$-continuous, if for every net $\{x_\alpha\}$ in $W_1$ with $x_\alpha \xrightarrow{\tilde{o}}0$ it follows that $T(x_\alpha)\xrightarrow{\tilde{o}}0$ holds in $W_2$.
		\item $\sigma$-$\tilde{o}$-continuous,
		 if for every sequence $\{x_n\}$ in $W_1$ with $x_n \xrightarrow{\tilde{o}}0$ it follows that $T(x_n)\xrightarrow{\tilde{o}}0$ holds in $W_2$.
	\end{enumerate}
We show the collection of all $\tilde{o}$-continuous ($\sigma$-$\tilde{o}$-continuous) operators between two ordered vector spaces $W_1,W_2$, by $L_{\tilde{o}c}(W_1,W_2)$ ($L_{\sigma \tilde{o}c}(W_1,W_2)$).
\end{definition}
Example 1.8 of \cite{11} shows that the classes of $\tilde o$-continuous operators and order continuous operators between two ordered vector spaces $W_1$ and $W_2$ are different.
\begin{example}\label{jhg}
		Let $B$ be a projection band of ordered vector space $W$ and $ P_{B}$
	the corresponding band projection.  Let $\{x_\alpha\}\subseteq W$ and $x_\alpha \xrightarrow{\tilde o}0$ in $W$. There exists a net $\{y_\beta\}\subseteq W$ that $y_\beta \downarrow 0$ in $W$. For each $\beta$ there is $\alpha_0$ such that $P_B x_\alpha \leq y_\beta$ for each $\alpha\geq \alpha_0$. It is clear that $y_\beta \xrightarrow{o}0$ in $W$ and there exists a net $\{z_\beta\}\subseteq B$ that $P_B(y_\beta)\leq z_\beta$ for each $\beta$. We have $P_B(x_\alpha) = P_B(P_B(x_\alpha))\leq P_B(y_\beta)\leq z_\beta$ whenever $\alpha \geq \alpha_0$.

\end{example}
\begin{theorem}\label{uyi}
	Let $T:W_1\rightarrow W_2$ be an operator between two ordered vector spaces.
	\begin{enumerate}
		\item If $0\leq T$ is order continuous, then $T$ is $\tilde{o}$-continuous.\label{ali}
		\item If $W_2$ is a Dedekind complete vector lattice and $T$ is $\tilde o$-continuous, then $T$ is order continuous.\label{alib}
		\item 	If $W_1$ is directed with RDP and $W_2$ is a Dedekind complete vector lattice. Then operator $T: W_1 \rightarrow W_2$ is order continuous iff $T$ is $\tilde o$-continuous.
	\end{enumerate} 
\end{theorem}
\begin{proof}
	\begin{enumerate}
\item	Let $\{x_\alpha\}\subseteq W_1$ be a net such that $x_\alpha \xrightarrow{\tilde{o}}0$. There exists a net $\{y_\beta\}$ in $W_1$ such that $y_\beta \downarrow 0$ and for every $\beta$ there exists $\alpha_0$ such that $\pm x_\alpha \leq y_\beta$ whenever $\alpha \geq \alpha_0$. Due to $T$ being positive we obtain $\pm T(x_\alpha)\leq T(y_\beta)$. Since $T$ is positive and order continuous, hence by Lemma 7 of \cite{2}, $Ty_\beta \downarrow 0$ in $W_2$. It follows $T(x_\alpha)\xrightarrow{\tilde{o}}0$.
\item Let $\{x_\alpha\}\subseteq W_1$ be a net such that $x_\alpha \xrightarrow{o}0$ in $W_1$. It is clear that $x_\alpha \xrightarrow{\tilde o}0$ in $W_1$. By assumption $Tx_\alpha \xrightarrow{\tilde o}0$ in $W_2$. Since $W_2$ is a Dedekind complete vector lattice, $Tx_\alpha \xrightarrow{o}0$ in $W_2$ (see page 288 of \cite{11}). Hence $T$ is order continuous.
\item Let $T\in L_{oc}(W_1,W_2)$.	By Proposition 11 of \cite{2} we have $L_{oc}^\diamond (W_1,W_2) = L_{oc}(W_1,W_2)$. Therefore $T = T_1 - T_2 $ such that $T_1$ and $T_2$ are positive and order continuous. By \ref{ali}, $T_1 , T_2$ are $\tilde o$-continuous and therefore $T$ is $\tilde o$-continuous.\\
Conversely, it is clear by  \ref{alib}.
\end{enumerate}
\end{proof}
In this section we want to bring up two propositions similar to Theorem 13 and Proposition 20 of \cite{2}, respectively.
\begin{proposition}
	Let $W_1$ be a directed ordered vector space and $W_2$ be a pre-Riesz space with a vector lattice cover $(W_3,i)$. If $T\in L_{\tilde{o} c}(W_1,W_2)$, then $ioT \in L_{\tilde{o} c}(W_1,W_3)$.
\end{proposition}
\begin{proof}
	The proof has an argument similar to Theorem 13 of \cite{2}.
\end{proof}
\begin{proposition}
	Let $V_1$ and $V_2$ be two pre-Riesz spaces and $(W_1 , i_1)$, $(W_2 , i_2)$ be their vector lattice covers, respectively. Let a positive  operator $T: V_1 \rightarrow V_2$ has a positive linear extension $S: W_1 \rightarrow W_2$, i.e. $ Soi_1 = i_2 o T$. If $S\in L_{\tilde oc}(W_1,W_2)$, then $T\in L_{{\tilde oc}}
	(V_1,V_2)$.
\end{proposition}
\begin{proof}
	Let $\{x_\alpha\}$ be a net in $V_1$ with $x_\alpha \downarrow 0 $ in $V_1$. By Lemma 1(i) of \cite{2}, the infimum of the set $\{i_1 (x_\alpha): \alpha \in A\}$ exists in $W_1$ and equals $0$. It is clear that $i_1(x_\alpha)\xrightarrow{\tilde o}0$ in $W_1$. As $S$ is ${\tilde o}$-continuous, $S(i_1(x_\alpha))\xrightarrow{\tilde o} 0$. Therefore there exists a net $\{y_\beta\}\subseteq W_2$ such that $y_\beta \downarrow 0$ and for each $\beta$ there is an $\alpha_0$ such that $\pm S(i_1(x_\alpha)) \leq y_\beta$ whenever $\alpha \geq \alpha_0$. Clearly $i_2 (T(x_\alpha))\downarrow$. We have $0 \leq i_{2}(T(x_\alpha)) = S(i_{1}(x_\alpha))\leq y_\beta$ in $W_2$.  It is clear that $i_2 (T(x_\alpha))\downarrow 0$ in $W_2$ and by Lemma 1(ii) of \cite{2}, $T(x_\alpha)\downarrow 0$ in $V_2$.	By Lemma 7 of \cite{2}, $T$ is order continuous and therefore by Theorem \ref{uyi}, $T$ is $\tilde o$-continuous.
\end{proof}

\section{Order convergence in semi-order vector spaces }

A vector space $V$ is called semi-order vector space (in short semi-order space), if there exist an ordered vector space $W$ and an operator $T$ from $V$ into $W$. In this way, we say that $V$ is semi-order space with respect to $\{W, T\}$. 
	 A net $\{x_\alpha\}\subseteq V$ is said to be ${\{W,T\}}$-order convergent (resp. $T(V)$-order convergent)  to a vector $x\in V$ in short we write $x_\alpha\xrightarrow {\{W, T\}}x$ (resp. $x_\alpha\xrightarrow {T(V)}x $), whenever there exists a net $\{y_\beta\}$ in $W$ (resp. $T(V)$)  satisfying $y_\beta \downarrow 0$ in $W$ (resp. $T(V)$)  and for each $\beta$, there exists $\alpha_0$ such that $\pm (Tx_\alpha -Tx) \leq y_\beta$ whenever $\alpha \geq \alpha_0$.\\
	  Whenever $V$ is a subspace of $W$ and $T$ is inclusion map, we use symbol $x_\alpha \xrightarrow{Wo}x$ in $V$ instead of symbol $x_\alpha  \xrightarrow{\{W,T\}}x$ in $V$.\\  
In this section, we study some properties of  $\{W,T\}$-order convergent nets.

 Let $K^\prime$ be a cone in $W$. Obviously,  $K^\prime \cap T(V)= K^{\prime\prime}$ is a cone of $T(V)$. Then there exists $K\subseteq V$ with $T(K)=K^{\prime\prime}$. Now if $\ker T = \{0\}$,  then  $K$ is a cone of $V$. It means that if $V$ is a semi-order vector space with respect
  to $\{W,T\}$ that $\ker T = \{0\}$, then $V$ is an ordered vector space. \\

 We say that $V$ has order properties, when $T(V)$ has these order properties, for example  see the following definition for some of them.
 
\begin{definition}\label{puy}
Assume that $V$ is a semi-order space with respect to $\{W,T\}$.
\begin{enumerate}
\item For each $x\in V$, we define $x\geqslant_V 0$ whenever $Tx\geqslant 0$ ($\geqslant_V$ is named semi-order in $V$).
\item A subset $M$ of $V$ is $\{W,T\}$-order closed, $\{W,T\}$-order bounded in  $V$ whenever $T(M)$ is    order closed, order bounded in  $W$, respectively.
\item A subspace $B\subseteq V$ is $\{W,T\}$-order dense, $\{W,T\}$-ideal, $\{W,T\}$-band in  $V$ whenever $T(B)$ is   order dense, ideal, band in  $W$, respectively. The operator $P_B : V\rightarrow B$ defined via $P_B(x)=x_1$ where $P_{TB}(Tx)=Tx_1$ that $x_1 \in B$ and $P_{TB}: W\rightarrow TB$ is a band projection, is a band projection on $V$.
\item Let  a net $\{x_\alpha\}\subseteq V$. $x_\alpha\downarrow x$ whenever $Tx_\alpha\downarrow Tx$ in $W$. 
\item A sequence $\{x_n\}\subseteq V$ is said to be  $\{W,T\}$-disjoint in $V$, if $\{Tx_n\}$ is disjoint  sequence in $W$.
\end{enumerate}
\end{definition}
For a semi-order space $V$, order convergence is depended to ordered vector space $W$ and  operator $T$.
  The following example, shows that ${\{W,T\}}$-order convergence depends on both of $W$ and $T$.
 
\begin{example}\label{elin}
		Consider the standard basis $\{e_n\}$ of $c_0$. $c_0$ is a subspace of $\ell^\infty$ and $e_n\xrightarrow{\{\ell^\infty o\} }0$ in $c_0$, but $\{e_n\}$ is not $I(c_0)$-order convergent to $0$ in $c_0$. Now set an operator $T$ from $c_0$ into $\ell^\infty$ with $Te_n=a_n$ whenever $a_n=(n,n,n,...n, 0,0,0,...)$,  the first $n$ places are occupied with $n$ and the rest is zero. It is clear that $\{e_n\}$ is not ${\{\ell^\infty,T\}}$-order convergent to zero.
\end{example}
\begin{lemma}
	Assume that  $V$ is a semi-order space  with respect to $\{W,T\}$. Then we have  the following assertions.
\end{lemma}
\begin{enumerate}
	\item $x_\alpha \xrightarrow{\{W, T\}}x$ in $V$ iff $(x_\alpha - x) \xrightarrow{\{W, T\}}0$ in $V$.
\item If $0\leq_V x_\alpha \xrightarrow{\{W, T\}}x$ in $V$, then $0 \leq_V x$.
\item If for each $\alpha$, $ x_\alpha \leq_V y$ in $V$ and $x_\alpha \xrightarrow{\{W, T\}}x$ in $V$, then $ x \leq_V y$.
\item If $x_\alpha\xrightarrow{\{W, T\}}x$ and $x_\alpha \xrightarrow{\{W, T\}}y$ in $V$, then $x = y$.
\item If $x_\alpha \xrightarrow{\{W, T\}}x$ and $y_\alpha \xrightarrow{\{W, T\}}y$ in $V$, then $ \lambda x_\alpha + \mu y_\alpha \xrightarrow{\{W, T\}} \lambda x + \mu y$ in $V$ for all $ \lambda , \mu \in\mathbb{R}$.
\item If $x_\alpha \xrightarrow{\{W, T\}}x$, $z_\alpha \xrightarrow{\{W, T\}}z$ and $ x_\alpha \leq_V z_\alpha$ for all $\alpha$ in $V$, then $ x \leq_V z$.
\end{enumerate}
\begin{proof}
By using of Definition \ref{puy} and the proof of  Lemma \ref{eli}, the proof is complete.
\end{proof}
\begin{theorem}
	\begin{enumerate}
\item	Let  $W$ be an order dense subspace of ordered vector space $U$. If $\{x_\alpha\}\subseteq V$ and $x_\alpha \xrightarrow{\{W, T\}}0$ in $V$, then $x_\alpha \xrightarrow{\{U, T\}}0$ in $V$.\label{rew}
\item If $V$ is order dense in $W$ and $x_\alpha \xrightarrow{T(V)}x$ in $V$, then $x_\alpha \xrightarrow{\{W,T\}}x$ in $V$.
\item Assume that  $V$  is semi-order space with respect to both $\{{{W_1}}, T\}$ and $\{W_2, T\}$ such that  ${{W_1}}$  is an ideal of vector lattice $W_2$. If $\{x_\alpha\}$ is $\{W_1,T\}$-order bounded in ${{V}}$, then $x_\alpha \xrightarrow{\{W_2, T\}}0$ in $V$ implies $x_\alpha\xrightarrow{\{{{W_1}}, T\}}0$ in $V$.
	\item Let $W$ be a vector lattice, $I$ be a $\{W,T\}$-ideal in $V$ and $\{x_\alpha\} \subseteq I$. $x_\alpha \xrightarrow{\{W,T\}}x$ in $V$ iff $x_\alpha \xrightarrow{\{W,T\mid_I\}}x$ in $I$.
	\item Let $(U,i)$ be a vector lattice cover of pre-Riesz space $W$ and $\{x_\alpha\}\subseteq V$. Then $x_\alpha \xrightarrow{{\{W,T\}}}0$ in $V$ implies $x_\alpha \xrightarrow{\{U, ioT\}}0$ in $V$.
\end{enumerate}
	\end{theorem}
\begin{proof}
\begin{enumerate}
\item  Let $\{x_\alpha\}\subseteq V$ and $x_\alpha \xrightarrow{\{W,T\}}x$ in $V$, therefore there exists a net $\{y_\beta\}\subseteq W$ that $y_\beta \downarrow 0$ in $W$. By Proposition 5.1 of \cite{3}, $y_\beta \downarrow 0$ in $U$. Hence the proof is complete.
\item By assumption and by Definition \ref{puy}, $T(V)$ is order dense in $W$. Therefore the proof is clear by \ref{rew}.
\item  Assume that   $\{x_\alpha\}\subseteq V$  and $x_\alpha \xrightarrow{\{W_2, T\}}0$ in $V$.
 Then there exists $\{y_\beta\}\subseteq W_2$ satisfying $y_\beta \downarrow 0$ and for each $\beta$ there exists $\alpha_0$ such that $ |Tx_\alpha|\leq y_\beta$ whenever $\alpha \geq \alpha_0$. By assumotion and by Definition \ref{puy}, there exists a $u\in {{W_1}}^+$ such that $|Tx_\alpha| \leq u$. Since  ${{W_1}}$ is an ideal of $W_2$, $\{u \wedge y_\beta\}\subseteq {{W_1}}$.  It is clear that  $u \wedge y_\beta \downarrow 0$  in ${{W_1}}$. It is obvious   that for each $\beta$ there exists $\alpha_0$ that $|Tx_\alpha | \leq u \wedge y_\beta$ whenever $\alpha \geq \alpha_0$. It follows that $x_\alpha \xrightarrow{\{{{W_1}}, T\}}0$ in $V$. 
 \item  Let $\{x_\alpha\} \subseteq I$ and $x_\alpha \xrightarrow{\{W,T\}}x$ in $V$, that means that $T(x_\alpha)\xrightarrow{\tilde o}T(x)$ in $W$. By Definition \ref{puy}, $T(I)$ is an ideal in $W$. By Lemma 2.5 of \cite{4}, $T(x_\alpha)\xrightarrow{\tilde o}T(x)$ in $T(I)$. Hence $x\in I$.\\
 Conversely, it is clear that if $x_\alpha \xrightarrow{\{W,T\mid_I\}}x$ in $I$, then $x_\alpha \xrightarrow{\{W,T\}}x$ in $V$.
 \item 	Let $\{x_\alpha\}\subseteq V$ and $x_\alpha \xrightarrow{{\{W,T\}}}0$ in $V$. Then there exists a net $\{y_\beta\}\subseteq W$ such that $y_\beta \downarrow 0$ in $W$ and for each $\beta$ there exists $\alpha_0$ that $\pm (T(x_\alpha - x))\leq y_\beta$ whenever $\alpha \geq \alpha_0$. Since $W$ is order dense in $U$, therefore by Proposition 5.1 of \cite{3}, $y_\beta \downarrow  0$ in $U$. By Lemma 1 of \cite{2}, $i(y_\beta)\downarrow 0$ in $U$. Because $i$ is positive, we have $\pm i(T(x_\alpha - x))\leq i(y_\beta)$.
\end{enumerate}	
\end{proof} 




Assume that $V$ is a semi-order space with respect  to $\{W,T\}$. A set  $A\subseteq V$ is $\{W, T\}$-closed (resp. $T(V)$-closed) in $V$, if for any net $\{x_\alpha\}\subseteq A$ and $x\in V$ with $x_\alpha\xrightarrow{\{W, T\}}x$ (resp. $x_\alpha\xrightarrow{T(V)}x$) in $V$, one has $x\in A$.
\begin{proposition}
Let $T$ be onto, one-to-one and $A\subseteq V$. Then $A$ is $\{W,T\}$-order closed iff $T(A)$ is $\tilde{o}$-closed in $W$
\end{proposition}
\begin{proof}
	Let $\{x_\alpha\} $ be a net in $A$ and $x_\alpha \xrightarrow{\{W,T\}}x$ in $V$. It follows that  there exists a net $\{y_\beta\}\subseteq W$ such that $y_\beta \downarrow 0$ and for each $\beta$ there exists $\alpha_0$ that $\pm (Tx_\alpha - Tx) \leq y_\beta$ whenever $\alpha \geq \alpha_0$. Since $\{Tx_\alpha\}\subseteq T(A)$ and $T(A)$ is $\tilde o$-closed in $W$, therefore $Tx\in  T(A)$. It follows that $x\in A$.\\
 Conversely, let $\{T(x_\alpha)\}\subseteq T(A)$ and $Tx_\alpha\xrightarrow{\tilde o}y$. Since $T$ is onto, there exists $x\in V$ that $Tx = y$.  It is clear that $\{x_\alpha\}\subseteq A$, $x_\alpha \xrightarrow{\{W,T\}}x$ in $V$. By assumption $x\in A$. So $Tx \in T(A)$. 
\end{proof}
\begin{proposition}
	Let $V$ be a semi-order space with respect to $\{W,T\}$ where $W$ is a vector lattice. Each $\{W,T\}$-disjoint and $\{W,T\}$-order bounded sequence $\{x_n\}\subseteq V$ is $\{W,T\}$-order null.
\end{proposition}
\begin{proof}
By Definition \ref{puy}, $\{Tx_n\}$ is order bounded and disjoint in $W$. Therefore by Corollary 3.6 of \cite{4}, $Tx_n \xrightarrow{uo}0$ in $W$. Since $\{Tx_n\}$ is order bounded in $W$, we have $Tx_n \xrightarrow{\tilde o}0$ in $W$. Hence there exists a sequence $\{y_m\}\subseteq W$ such that $y_m \downarrow 0$ in $W$ and for every $m$ there exists $n_0$ such that $\pm Tx_n \leq |Tx_n|\leq y_m$  whenever $n \geq n_0$. Thus $x_n \xrightarrow{\{W,T\}}0$ in $V$.
\end{proof}

\section {Semi-order continuous operators}

Assume that $V_1$ and $V_2$ are semi-order spaces with respect to 	 $\{{{W_1}}, T_1\}$ and $\{W_2, T_2\}$, respectively. 
An operator $S$ from $V_1$ into $V_2$ is positive, when $x>_{V_1}0$ implies $Sx>_{V_2}0$

\begin{definition}
Assume that $V_1$ and $V_2$ are semi-order spaces with respect to 	 $\{{{W_1}}, T_1\}$ and $\{W_2, T_2\}$, respectively. An operator $S$ from $V_1$ into $V_2$ is called
\begin{enumerate}
	\item  semi-order continuous, if  $x_\alpha\xrightarrow {\{{{W_1}}, T_1\}}x$ implies $Sx_\alpha\xrightarrow {\{W_2, T_2\}}Sx$ whenever $\{x_\alpha\}\subseteq V_1$.
	\item  $\sigma$-semi-order continuous, if  $x_n\xrightarrow {\{{{W_1}}, T_1\}}x$ implies $Sx_n\xrightarrow {\{W_2, T_2\}}Sx$ whenever $\{x_n\}\subseteq V_1$.
\end{enumerate}
\end{definition}
It is obvious that an operator $S$ from $V_1$ into $V_2$ is called semi-order continuous if and only if $x_\alpha\xrightarrow {\{{{W_1}}, T_1\}}x$ implies $x_\alpha\xrightarrow {\{W_2, T_2S\}}x$ whenever $\{x_\alpha\}\subseteq V_1$.
Since semi-order continuity of an operator $S$ between two  semi-order spaces depends on ${{W_1}}$ and $W_2$, the collection of all semi-order continuous operators between two semi-order spaces $V_1$ and $V_2$ will be denoted by  $L_{o_{{{W_1}}W_2}}(V_1, V_2)$. Similarly, $ L_{\sigma o_{{{W_1}}W_2}}(V_1, V_2)$ will denote the collection of all operators from $V_1$ to $V_2$ that are $\sigma$-semi-order continuous. Whenever $S\in L_{o_{{{W_1}}W_2}}(V_1, V_2)$ (resp. $ L_{\sigma o_{{{W_1}}W_2}}(V_1, V_2)$), we say that, $S$ is $o_{W_1W_2}$-continuous (resp, $\sigma$-$o_{W_1W_2}$-continuous).\\
Here are some examples of semi-order continuous operators.


\begin{example}
\begin{enumerate}
	\item  Let $V$ be a pre-Riesz space and $W$ is its vector lattice cover that $W$ is Archimedean. The inclusion map $I : V \rightarrow W$ is $o_{V^\delta W^\delta}$-continuous ($V^\delta$, $W^\delta$ are Dedekind completions of $V$,$W$ respectively). Let $\{x_\alpha\}\subseteq V$ and $x_\alpha\xrightarrow{\{V^\delta o\}}0$ in $V$. Then by Theorem 2 of \cite{2}, $Ix_\alpha \xrightarrow{\{W^\delta o\}}0$ in $W$.
	\item Let $V$ be a semi-order space with respect to $\{W,T\}$, $B$ be a band of $V$ and $P_B:V\rightarrow B$ be a band projection. If $\{x_\alpha\}\subseteq V$ and $x_\alpha \xrightarrow{\{W,T\}}0 $ in $V$, then $Tx_\alpha\xrightarrow{\tilde o}0$ in $W$. By Definition \ref{puy}, $T(B)$ is a band in $W$. We consider the band projection $P_{TB}:W\rightarrow TB$. By Example \ref{jhg}, $P_{TB}(Tx_\alpha)\xrightarrow{\tilde o}0$ in $TB$. It is clear that $T|_B P_B = P_{TB}T$. Therefore $P_B(x_\alpha)\xrightarrow{\{T|_B(B)\}}0$ in $B$. So $P_B$ is $o_{W,TB}$-continuous.
	\item Let $f$ be a continuous functional on Banach lattice $E$. Since $\mathbb{R}$ is a $KB$-space so by Theorem 4.60 of \cite{1}, $c_0$ is not embeddable in $E$. By Theorem 4.63 of \cite{1}, there exist a $KB$-space $F$, lattice homomorphism $Q:E\rightarrow F$ and functional $g:F\rightarrow \mathbb{R}$ such that $ f = g o Q$. If  $I$ is identity map on $\mathbb{R}$, we have $Iof = g o Q$. If $\{x_\alpha\}\subseteq E$ and $x_\alpha\xrightarrow{\{F,Q\}}0$ in $E$, then $Q(x_\alpha)\xrightarrow{\tilde o}0$ in $F$. Since $F$ has order continuous norm, therefore $Q(x_\alpha)\xrightarrow{\|.\|}0$ in $F$. So $g(Q(x_\alpha))\xrightarrow{\|.\|}0$ and hence $g(Q(x_\alpha))\xrightarrow{\tilde o}0$ in $\mathbb{R}$. Therefore we have $f(x_\alpha)\xrightarrow{\{\mathbb{R}o\}}0$. Hence $f$ is  $o_{F,\mathbb{R}}$-continous.
\end{enumerate}
\end{example}
\begin{theorem}\label{ebadi}
	Let $V_1 , V_2$ be two pre-Riesz spaces, $(W_1,i_1), (W_2,i_2)$ be their vector lattice covers, respectively and $S:V_1 \rightarrow V_2$ be a positive operator. 
	\begin{enumerate}
		\item If $S$ has extension positive and order continuous $\tilde{S} : W_1\rightarrow W_2$ that $i_2oS = \tilde{S}o i_1$, then $S\in L_{o_{{{W_1}}W_2}}(V_1, V_2)$.
		\item 
		$ S\in L_{o_{{{W_1}}W_2}}(V_1, V_2)$ iff  $ x_\alpha \downarrow 0$ in ${{V_1}}$ implies $ Sx_\alpha \downarrow 0$ in $V_2$ for each net $\{x_\alpha\}\subseteq V_1$.
		\end{enumerate}
\end{theorem}
\begin{proof}
	\begin{enumerate}
	\item	 Let $\{x_\alpha\}\subseteq V_1$ and $x_\alpha \xrightarrow{\{W_1 , i_1\}}0$ in $V_1$. Then there exists a net $\{y_\beta\}\subseteq W_1$ such that $y_\beta \downarrow 0$ and for each $\beta$ there exists $\alpha_0$ that $\pm i_1(x_\alpha)\leq y_\beta$ whenever $\alpha \geq \alpha_0$. Since $\tilde{S}$ is order continuous and positive, therefore by Theorem 1.56 of \cite{1}, $\tilde{S}(y_\beta) \downarrow 0$ in $W_2$ and it is clear that $\tilde{S}(\pm i_1(x_\alpha))\leq 
	\tilde{S}(y_\beta)$. Hence by assumption $S(x_\alpha)\xrightarrow{\{W_2,i_2\}}0$.
	\item 
	Suppose that $0\leq S\in L_{o_{{{W_1}}W_2}}(V_1, V_2)$ and $\{x_\alpha\}\subseteq V_1$ with $x_\alpha \downarrow 0$ in ${{V_1}}$. Then by Lemma 1 of \cite{2}, $i_1(x_\alpha)\downarrow 0$ in $W_1$. It is clear that
	$x_\alpha \xrightarrow{\{{{W_1 ,i_1}}\}}0$ in $V_1$.
	By assumption, we have    $Sx_\alpha \xrightarrow{\{W_2,i_2\}}0$ in $V_2$. It follows that there exists a net $\{y_\beta\} \subseteq W_2$ satisfying, $y_\beta \downarrow 0$ in $W_2$ and for each $\beta$ there exists $\alpha_0$ such that $\pm i_2Sx_\alpha \leq y_\beta$ whenever $\alpha \geq \alpha_0$. It follows that $ i_2Sx_\alpha \leq \inf y_\beta =0$ for each $\alpha$ and so $ i_2Sx_\alpha \downarrow 0$ in $W_2$ and therefore by Lemma 1 of \cite{2}, $Sx_\alpha \downarrow 0$ in $V_2$.\\
	Conversely, let $\{x_\alpha\}\subseteq V_1$ and $x_\alpha \xrightarrow{\{W_1,i_1\}}0$ in $V_1$. There exists a net $\{y_\beta\}\subseteq W_1$ such that $y_\beta \downarrow 0$ in $W_1$ and for each $\beta$ there exists $\alpha_0$ such that $\pm i_1(x_\alpha)\leq y_\beta$ whenever $\alpha \geq \alpha_0$. Hence $i_1 x_\alpha \downarrow 0$ in $W_1$. By Lemma 1 of \cite{2}, $x_\alpha \downarrow 0$ in $V_1$. By assumption $Sx_\alpha\downarrow 0$ in $V_2$ and by Lemma 1 of \cite{2}, $i_2 S(x_\alpha)\downarrow 0$ in $W_2$. Therefore $S(x_\alpha)\xrightarrow{\{W_2,i_2\}}0$ in $V_2$.
	\end{enumerate}
\end{proof}





Assume that $V_1$ and $V_2$ are semi-order spaces with respect to 	 $\{{{W_1}}, T_1\}$ and $\{W_2, T_2\}$, respectively. 
An operator $S:V_1\rightarrow V_2$ is semi-order bounded, if $S(A)$ is  $\{W_2,T_2\}$-bounded for each  $\{W_1,T_1\}$-bounded set $A\subseteq V_1$.

 If  $S:V_1\rightarrow V_2$  semi-order bounded, then we write that $S$ is  ${{W_1}}W_2$-order bounded (for short, $o_{{{W_1}}W_2}$-bounded). 
 \\
We denote the collection of all  $o_{{{W_1}}W_2}$-bounded operators $S:V_1\rightarrow V_2$ by $L_{b_{{{W_1}}W_2}}(V_1, V_2)$.

\begin{theorem}\label{kat}
	Let $V_1, V_2$ be  subspaces of ordered vector spaces ${{W_1}}, W_2$, respectively, such that $W_2$ is  Archimedean Dedekind complete vector lattice. Then, $L_{o_{{{W_1}}W_2}}(V_1, V_2)$ is subspace of $ L_{b_{{{W_1}}W_2}}(V_1, V_2)$.
\end{theorem}
\begin{proof}
 Let $S:V_1 \rightarrow V_2$ be an $o_{W_1 W_2}$-continuous operator. First we consider $y \in {{W_1}}_+$ and  $A = V_1\cap [0,y]$.
 Let $I= \mathbb{N}\times A$ be an index set with the lexicographical order. Namely, $(n,x)>(m,z)$ if and only if either one the following holds true.\\ 
$(1)$ $n>m$,\\
$(2)$ $n=m$ and $x>z$.\\
It is easy to check that $I$ is a directed set, so we may consider a net indexed by $I$. Let us set $\varphi_{(n,x)}=\frac{1}{n}x$ for all $x\in A$. Then we have $0\leq\varphi_{(n,x)}\leq \frac{1}{n}x$. It follows that $\varphi_{(n,x)}$ is $\{W_1\}$-order convergent to zero.
By assumption,  $S\varphi_{(n,x)}$ is $\{W_2\}$-order convergent to zero.
Then there exists a net $(y_\beta)_{\beta}$ such that 
$y_\beta\downarrow 0$ and for every $\beta$ there exists $(n,x)$ satisfying $ \pm S\varphi_{(m,z)}\leq y_\beta$ for all $(m,z)>(n,x)$.
Let us pick any $y_\beta$ and find corresponding index $(n,x)\in I$. Then, in particular, $\pm S\varphi_{(n+1,z)}\leq y_\beta$ for all $z\in A$. It follows that $-(n+1)y_\beta\leq Sz\leq (n+1)y_\beta$ for every 
$z\in A$. Thus, $S$ is $o_{{{W_1}}W_2}$-bounded operator from $V_1$ into $V_2$. 
\end{proof}

Let $V_1$ and $V_2$ be two semi-order spaces with respect to $\{W_1,T_1\}$ and $\{W_2,T_2\}$, respectively, where $W_1$ is directed. We define the directed part of $L_{o_{{{W_1}}W_2}}(V_1,V_2)$ by
$$L_{o_{{{W_1}}W_2}}^\diamond(V_1,V_2) := L_{o_{{{W_1}}W_2}}(V_1,V_2)_+ - L_{o_{{{W_1}}W_2}}(V_1,V_2)_+.$$
\begin{theorem}\label{po}
	Let $V_1$ and $V_2$ be two  pre-Riesz spaces with vector lattice covers $\{W_1,i_1\}$ and $\{W_2,i_2\}$,respectively, such that $V_1$ with a generating cone has the RDP and $V_2$, $W_2$ is Archimedean  Dedekind complete vector lattice. Then 
	\begin{enumerate}
		\item $L_{b_{{{W_1}}W_2}}(V_1,V_2)$ is a Dedekind complete vector lattice.\label{iu}
		\item $L_{o_{{{W_1}}W_2}}(V_1,V_2)$ is a vector lattice.
		\item  $L_{o_{{{W_1}}W_2}}(V_1,V_2)$ is a band in $L_{b_{{{W_1}}W_2}}(V_1,V_2).$
		\end{enumerate}
\end{theorem}
\begin{proof}
	\begin{enumerate}
		\item By Theorem 8 of \cite{2}, $L_b(V_1,V_2)$ is a Dedekind complete vector lattice. Since $V_1$ and $V_2$ are majorizing in ${{W_1}}$ and $W_2$ respectively, it is clear that $T\in L_{b_{{{W_1}}W_2}}(V_1,V_2)$ iff $T\in L_b(V_1,V_2)$. Therefore $L_{b_{{{W_1}}W_2}}(V_1,V_2)$ is a Dedekind complete vector lattice.
		\item Let $T\in L_{o_{{{W_1}}W_2}}(V_1,V_2)$. By Theorem \ref{kat}, $T\in L_{b_{{{W_1}}W_2}} (V_1,V_2)$. Therefore by \ref{iu}, $|T|$ exists and belongs to $L_{b_{{{W_1}}W_2}} (V_1,V_2)$. Let $\{x_\alpha\}\subseteq V_1$ and $x_\alpha \downarrow 0$. By Theorem \ref{ebadi}$(2)$, it is enough to show that $|T|x_\alpha\downarrow 0$ in $V_2$. Since $|T|$ is positive and $V_2$ is a Dedekind complete, so there exists a $z\in V_2$ that $|T|x_\alpha \downarrow z$ in $V_2$. Let $z\neq 0$. We have $i_2 |T|(x_\alpha)\downarrow i_2z$ in $W_2$. There exists a net $\{y_\beta\}\subseteq W_2$ that $y_\beta\downarrow 0$ and for each $\beta$ there exists $\alpha_0$ that $\pm i_2 (|T|x_\alpha -z)\leq y_\beta$ whenever $\alpha \geq \alpha_0$. We have $\pm i_2 (Tx_\alpha) \leq i_2 (|T|x_\alpha)\leq y_\beta \pm i_2 z$. So $\pm i_2 (Tx_\alpha - z)\leq y_\beta$ whenere $\alpha\geq \alpha_0$. This is a contradiction with $Tx_\alpha \xrightarrow{\{W_2,i_2\}}0$.  Therefore $|T|x_\alpha \downarrow 0$ in $V_2$.
\item  By Theorem \ref{kat}, $L_{o_{{{W_1}}W_2}}(V_1,V_2)$ is a subspace in  $L_{b_{{{W_1}}W_2}} (V_1,V_2)$. Let $T\in L_{o_{{{W_1}}W_2}}(V_1,V_2)$, $S\in L_{b_{{{W_1}}W_2}}(V_1,V_2)$ with $|S|\leq |T|$ and $\{x_\alpha\}\subseteq V_1$ with $x_\alpha \xrightarrow{\{W_1, i_1\}}0$ in $V_1$. We have $|T|x_\alpha \xrightarrow{\{W_2,i_2\}}0$.  With loss of generality, we can assume $0\leq x_\alpha$ for each $\alpha$. By enequality $\pm i_2 S(x_\alpha)\leq i_2 |S|x_\alpha \leq i_2 |T|x_\alpha$, $Sx_\alpha \xrightarrow{\{W_2,i_2\}}0$. So $L_{o_{{{W_1}}W_2}}(V_1,V_2)$ is an ideal in $L_{b_{{{W_1}}W_2}}(V_1,V_2)$. To see that the ideal $L_{o_{{{W_1}}W_2}}(V_1,V_2)$ is a band, the proof has similar argument of Theorem 1.57 \cite{1}.
\end{enumerate}\end{proof}

\end{document}